\documentclass{article}
\usepackage{graphicx} 
\usepackage{subcaption}
\usepackage{amsmath}
\usepackage{float} 
\usepackage{comment}
\usepackage[utf8]{inputenc}

\usepackage[a4paper, margin=1in]{geometry} 

\usepackage{geometry}
\usepackage{hyperref}
\usepackage[normalem]{ulem}
\usepackage{authblk}
\usepackage{lipsum}
\title{Social dynamics can delay or prevent climate tipping points by speeding the adoption of climate change mitigation}
\author[1,2]{Yazdan Babazadeh Maghsoodlo}
\author[2]{Madhur Anand}
\author[1]{Chris T. Bauch}

\affil[1]{Department of Applied Mathematics, University of Waterloo, Waterloo, Ontario, Canada,}
\affil[2]{School of Environmental Sciences, University of Guelph, Guelph, Ontario, Canada}

\begin{document}

\maketitle
\section{Abstract}


Social behaviour models are increasingly integrated into climate change studies, and the significance of climate tipping points for `runaway' climate change is well recognised. However, there has been insufficient focus on tipping points in social-climate dynamics. We developed a coupled social-climate model consisting of an Earth system model and a social behaviour model, both with tipping elements. The social model explores opinion formation by analysing social learning rates, the net cost of mitigation, and the strength of social norms. Our results indicate that the net cost of mitigation and social norms have minimal impact on tipping points when social norms are weak. As social norms strengthen, the climate tipping point can trigger a tipping element in the social model. However, faster social learning can delay or prevent the climate tipping point: sufficiently fast social learning means growing climate change mitigation can outpace the oncoming climate tipping point, despite social-climate feedback. By comparing high- and low-risk scenarios, we demonstrated high-risk scenarios increase the likelihood of tipping points. We also illustrate the role of a critical temperature anomaly in triggering tipping points. In conclusion, understanding social behaviour dynamics is vital for predicting climate tipping points and mitigating their impacts.


\section{Introduction}

Using various measurement techniques, we know with confidence that the Earth's average global temperature has increased by 1.1°C since 1800 due to anthropogenic activities \cite{zhou2021new}. For billions of years, Earth has developed intricate feedback mechanisms among various carbon reservoirs, balancing carbon emissions and absorption to achieve equilibrium. Examples of these feedback mechanisms include the ocean's capacity to absorb excess atmospheric carbon and forests on land capturing carbon from the atmosphere. \cite{williams2019carbon} \cite{heinze2019esd}. 
Since the beginning of industrialisation, excess human-caused carbon emissions have interfered with the function of these mechanisms \cite{majeed2020effects}. In addition to direct carbon emissions, human activities can have indirect effects that can exacerbate the situation. Recent studies have revealed a reduced capacity of oceans and land to absorb atmospheric carbon due to ocean acidification and forest fires. \cite{sivakumar2007climate} \cite{doney2009ocean}. This is among the causes of the rise in the global average temperature as more carbon is released into the atmosphere. Consequently, the ability of oceans and land to absorb carbon decreases, leaving more carbon in the atmosphere and forming a positive feedback loop. This positive feedback loop can increase the probability of irreversible events and a tipping point due to temperature rise \cite{marques2020climate}, \cite{eker2024cross}.

The primary way to combat global warming currently is mitigation through emissions reductions \cite{pathak2024climate}, \cite{fawzy2020strategies}, and based on the fraction of mitigators in a population, total carbon emissions can be controlled. In this view, Earth and people are seen as a coupled system, with people acting as elements of this generalised  system.\cite{liu2007complexity}, \cite{markkanen2019social}. This underscores the importance of behaviour models in understanding the mitigation process. Many behaviour models have been developed to provide a better understanding of how a population forms and exchanges opinions \cite{peralta2022opinion} \cite{acemoglu2011opinion} \cite{li2022game} \cite{zhu2022agent} \cite{biswas2011mean}.






In addition, in the past decades, there have been numerous efforts to develop a numerical structure for simulating climate models \cite{randall2007climate}. These models differ in the features and aspects that they include or exclude. Many of these models include behavioural models as well to offer a broad perspective of the future \cite{wood2011climate} \cite{10.1371/journal.pcbi.1007000}
\cite{10.1098/rspb.2021.1357}. Additionally, studying the possibility and existence of tipping points in the climate system has received high attention in previous years and many papers have explored different facets \cite{lenton2012arctic} \cite{lenton2008tipping} \cite{dakos2008slowing} \cite{bury2021deep}, including for climate tipping points \cite{wunderling2023global}.



 

Building on previous work \cite{beckage2022incorporating} \cite{bury2019charting} \cite{menard2021conflicts}, we aim to couple an earth system model with a behavioural model that has room for feedback mechanisms that can lead to tipping points. We aim to develop a model to understand the effects of events able to create positive feedback loops, like ocean acidification or forest fires, and the impact of different mitigation strategies on controlling the damage of these events by simulating the coupled model. This model is not intended to serve as a highly realistic forecast of future climate conditions. Instead, it is designed as a simplified tool to explore hypotheses, test assumptions, and investigate potential dynamic patterns. While it offers valuable insights, it does not aim to replace more comprehensive climate models.

We use the Earth System Model (ESM) and behavioural dynamics that were previously used to provide insight into the interaction between the climate system and social opinion formation \cite{lenton2000land} \cite{10.1371/journal.pcbi.1007000}. We introduce an additional term to simulate the effect of runaway carbon due to positive feedback mechanisms. Since this mechanism intensifies at higher temperatures and the probability of a tipping point increases in such situations, we construct this term to activate around a specific temperature threshold to simulate an irreversible tipping point. We also study the effect of social processes on the likelihood, strength, and timing of the climate tipping point. We have shown that the tipping point is preventable using specific mitigation strategies.

\section{Model Description}

\subsection{Behaviour Dynamics}
In this work, we use a behavioural dynamics model previously applied to social-climate systems \cite{10.1371/journal.pcbi.1007000}. Here, the population is assumed to be divided into two groups, mitigators and non-mitigators. Individuals in these two groups, decide whether to mitigate or not, based on the utility of each option. The utility of each option includes its associated cost, the costs associated with climate change, and the tendency to remain in the already chosen group. The fraction of the number of mitigators is represented by $x$ and the rate at which people can learn and switch their group is the social learning rate. The dynamic of $x$ is given by:
\begin{equation}
    \frac{dx}{dt} = \kappa x(1-x)[-\beta + f(T) + \delta(2x-1)]
\end{equation}

In this equation $\kappa$ is the social learning rate, $\beta$ is the net cost of mitigation, $\delta$ is the strength of social norms, and f(T) is the perceived costs associated with climate change, which is assumed to have the form of a sigmoid function:

\begin{equation}
    f(T) = \frac{f_{max}}{1+ e^{-\omega(T-T_{lim})}}
\end{equation}

The $f_{max}$ is the maximum cost, $\omega$ is the degree of non-linearity of the sigmoid function, and $T_{lim}$ is the critical temperature at which costs are most sensitive to change. Full details about deriving these equations are provided in the supplementary information.

\subsection{Earth system model}

To study the effect of social actions on climate change, We pair the social learning model with an earth system model \cite{lenton2000land}. Dynamics for atmospheric CO2 are modified to include an anthropogenic emission term that depends on the proportion of non-mitigators :

\begin{equation}
    \frac{dC_{at}}{dt} = \epsilon(t)(1-x) - P + R_{veg} + R_{so} - F_{oc} + \frac{R_{max}}{1+e^{R_0(T-T_c)}}
\end{equation}

$C_{at}$ is the deviation in the atmospheric $CO_2$ from pre-industrial levels, $\epsilon(t)$ is the baseline rate of $CO_2$ emission in the absence of mitigation. P is the rate of carbon uptake via photosynthesis, $R_{veg}$ is the outward carbon flux via plant respiration, $R_{so}$ is the outward carbon flux via soil respiration, and $F_{oc}$ is the net carbon uptake by the oceans. The last term is introduced to simulate and add the possibility of a tipping point (due to incidents like forest fires, and the ocean saturation) into our simulation and study. We assume this term to have a sigmoid function with $R_{max}$ to be the maximum value, $R_{0}$ be the degree of non-linearity of the sigmoid function, and $T_c$ be the critical temperature about which this term gets activated. The sigmoid function is chosen to add tipping point behaviour. In this paper, we use the \textit{modified model} for the results in which the tipping term was included in the simulations and \textit{baseline model} when it was off. Global surface temperature is assumed to evolve with the carbon cycle:
\begin{equation}
    c \frac{dT}{dt} = (F_d - \sigma T^4)a_E
\end{equation}

Where T is the deviation in global surface temperature from pre-industrial values and c is the thermal capacity specific heat capacity of the Earth’s surface. Also, $F_d$ is the net downward flux of radiation absorbed at the planet’s surface, $\sigma$ is the Stefan-Boltzmann constant, and $a_E$ is the Earth’s surface.

\subsection{Simulation}

In this work, we simulate our model from year 1800 to 2200. We kept the behaviour dynamic inactive from 1800 to 2017 and for $\epsilon$ in this period, we used the values of the historical record \cite{gilfillan2021cdiac}. For years after 2017, we projected the value of $\epsilon(t)$ using this equation:
\begin{equation}
    \epsilon(t) = \epsilon(2017) + \frac{(t-2017)\epsilon_{max}}{(t-2017+s)}
\end{equation}
Where $\epsilon_{max}$ is the saturation value and s is the half-saturation constant.
The initial values of all climate variables are zero since they represent the deviation from pre-industrial values. Social dynamics started in 2017 with a value of 5 percent for the initial fraction of mitigators ($x_0$).
The values of climate parameters are drawn from the original Earth system model \cite{lenton2000land}.
Social factors are more uncertain and have broader upper and lower limits. We applied the social parameter values from the reference paper to the behaviour model \cite{10.1371/journal.pcbi.1007000}.
For the tipping term, we varied $R_{max}$ across a broad range (0 to 5 GtC/year) to capture a wide spectrum of emission scenarios. For $T_c$ we used two values of 2 and 3 degrees Celsius to define high and low risk cases. For $R_0$ we used  5 $K^{-1}$ to ensure a sharp transition to the tipping point. The values of all the used parameters are provided in section 3 of the supplementary materials.

\section{Results}

\begin{figure}[t]
\includegraphics[width=\textwidth]{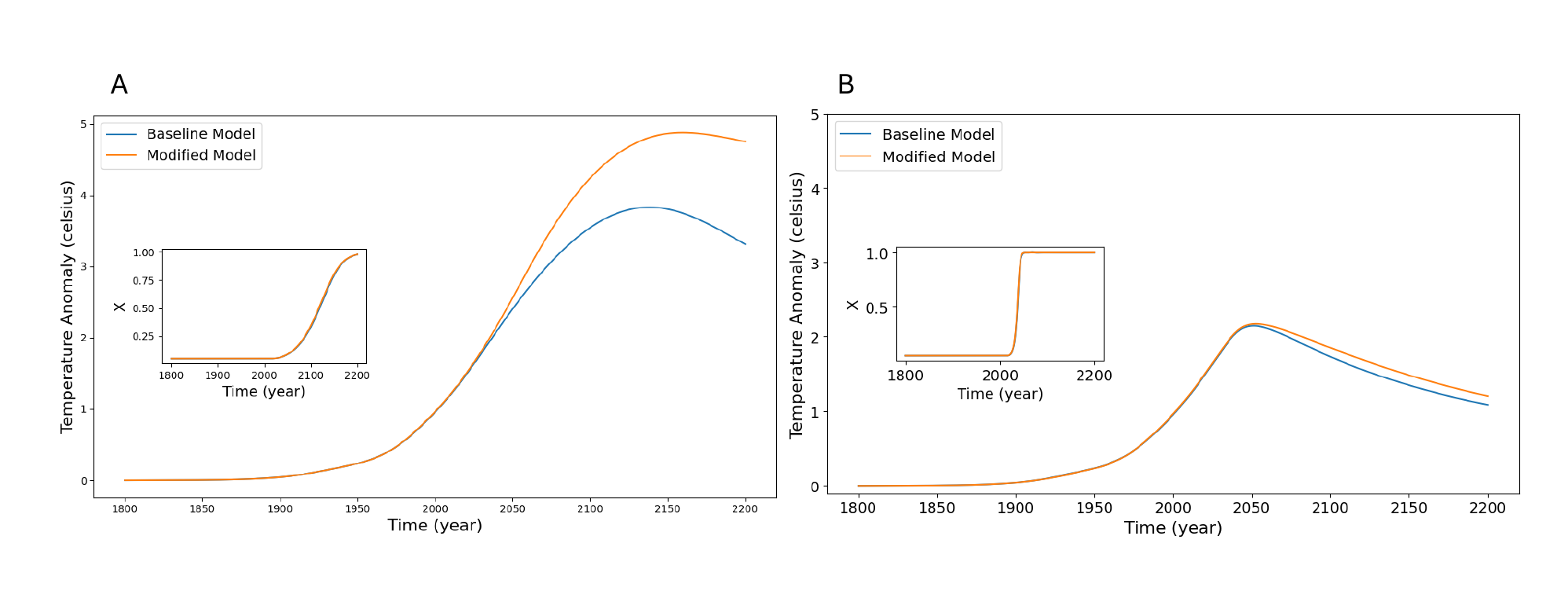}
\caption{\textbf{Possibility of avoiding climate tipping point through social actions}. In this figure, The Baseline model is illustrated in blue and the modified model with an additional tipping term is illustrated in orange. The dynamic of $x$ (fraction of mitigators) versus time is also depicted in both plots in a smaller box. Two plots differ in the values of $\kappa$ and $R_{max}$ (GtC/year) that were used: (0.01, 5) for plot A, where it shows tipping behaviour, and (0.1, 1) for plot B, where it does not.  Note the difference in vertical scale.}
\label{fig1}
\end{figure}

Simulating different scenarios, we found that it is possible to diffuse or delay the climate tipping point through social actions. In this model, $R_{max}$ (maximum value of the tipping term) represents the strength of the tipping term in the modified model, and $\kappa$ (social learning rate), $\beta$ (net cost of mitigation), and $\delta$ (strength of social norms) are parameters of the social behaviour model. Figure \ref{fig1} illustrates simulations of modified and baseline models for two examples of $\kappa$ and $R_{max}$. Comparing the two plots, we observe that the tipping term can lead to a significant difference between the modified and baseline models for certain choices of $\kappa$ and $R_{max}$, but not for others. This observation suggests that some social strategies can avoid the tipping point by deactivating the tipping term. To investigate and measure this further, we focus on the area under the curve (AUC) between the modified and baseline models as an indicator of the effect of the tipping term and how much difference it can make with different social variables. The difference in AUC represents the effect of the tipping term, which is the only difference between the modified and baseline models. Additionally, due to uncertainty about the value of critical temperature, we define two scenarios: 1) \textit{high risk}, in which the critical temperature is set at 2°C, and 2) \textit{low risk}, in which the critical temperature is set at 3°C.


\begin{figure}[!t]
\includegraphics[width=\textwidth]{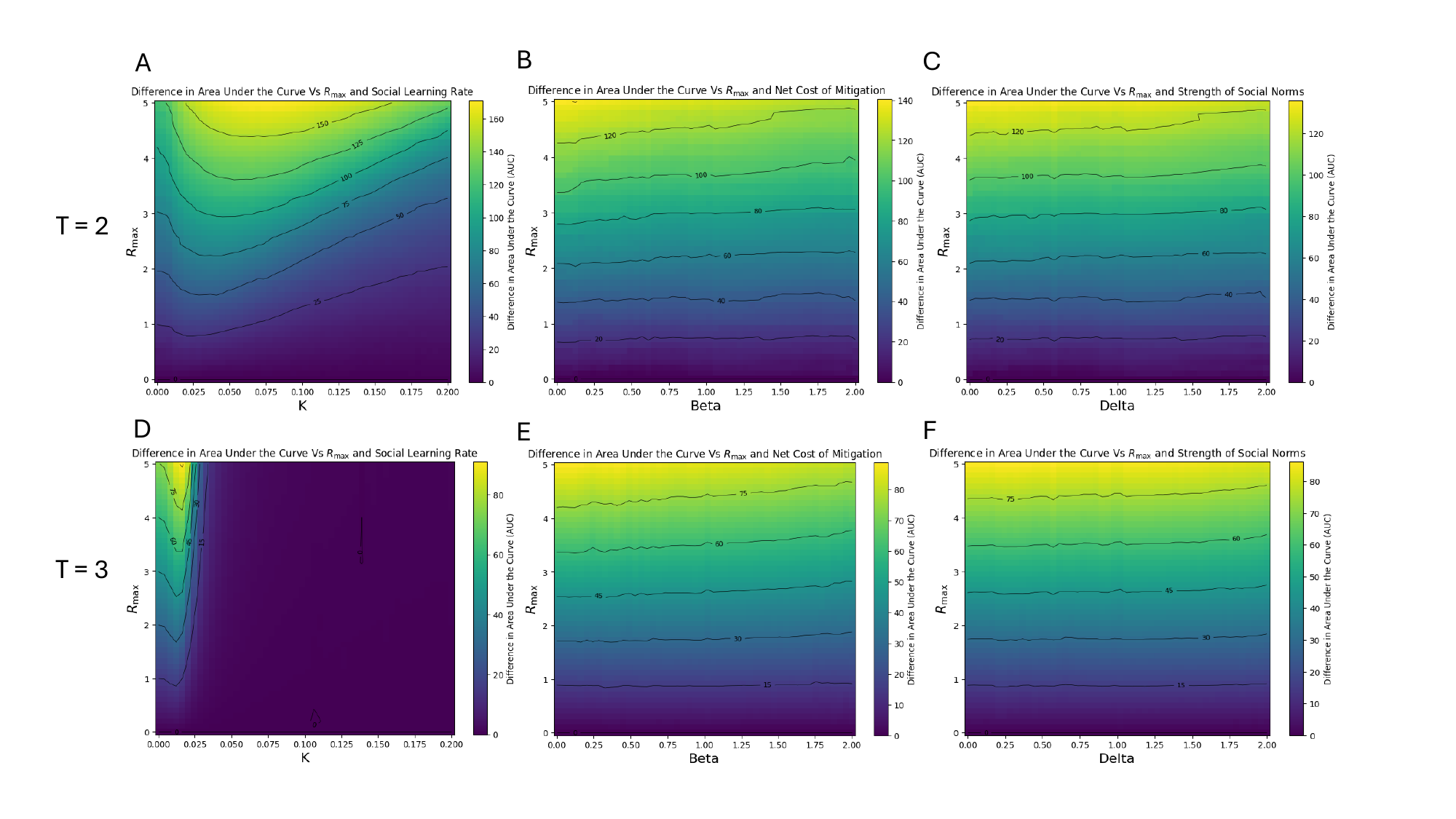}
\caption{\textbf{The difference in AUC versus social behaviour variables}. In this figure, the difference in AUC is depicted against social behaviour parameters and $R_{max}$ (GtC/year). The top plots represent the high-risk scenario ($T_c$ = 2), and the bottom ones correspond to the low-risk scenario ($T_c$ = 3). The plots illustrate the difference in AUC versus Social Learning Rate on the left, Net Cost of Mitigation in the middle, and Strength of Social Norm on the right.}
\label{fig2}
\end{figure}

Studying the difference in AUC versus social variables provides a more comprehensive insight into the potential to diffuse the tipping point. Figure \ref{fig2} illustrates this analysis, revealing that the difference in AUC does not rely heavily on $\beta$ (net cost of mitigation) and $\delta$ (strength of social norms). However, it depicts more complex behaviour when it comes to $\kappa$ (social learning rate). As shown in panels A and D, isoclines fit into V-shaped regions, meaning that for each value of $R_{max}$, there exist two values for $\kappa$ that result in the same difference in AUC, and there is one value of $\kappa$ between them that maximises it. The tipping term activates inside these V-shaped regions, suggesting the possibility of a tipping point. We can also compare the results in two predefined scenarios of high ($T_c = 2$) and low risk ($T_c = 3$) in Figure \ref{fig2}. As shown, panels A, B, and C represent the high-risk scenario, depicting more concerning results and a higher possibility of tipping points due to more elevated differences in AUC in panels B and C, as well as larger V-shaped regions. On the other hand, panels D, E, and F show that in low-risk scenarios, the V-shaped regions are smaller, and the difference in AUC is less pronounced. Since $\kappa$ has the greatest influence on the tipping term, we'll now focus primarily on results as a function of $\kappa$.


\begin{figure}[t]
\hspace*{-0.7cm}
\includegraphics[width=1\textwidth]{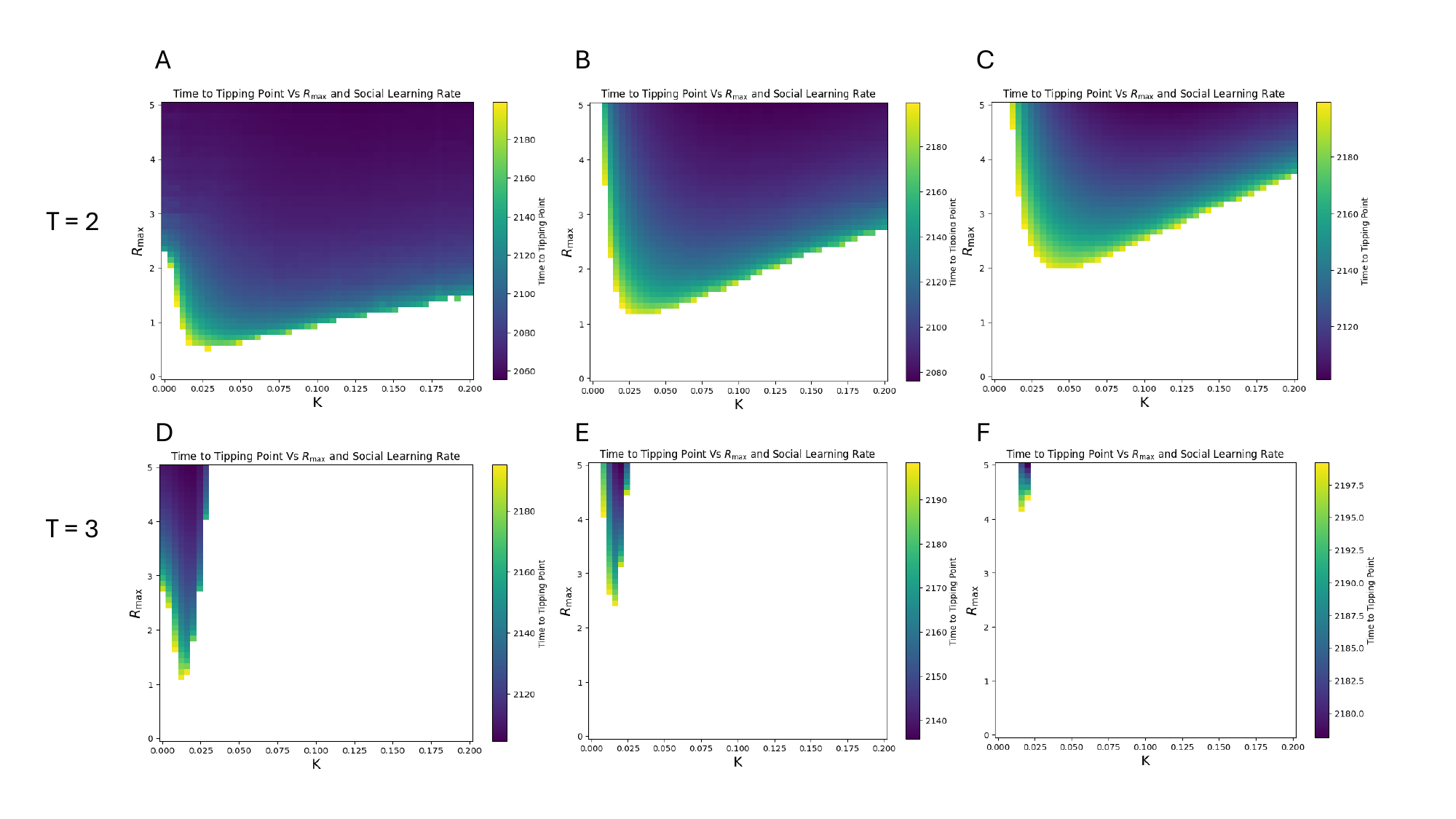}
\caption{\textbf{Time to the tipping point for different emission and mitigation scenarios}. Time to the tipping point is plotted against $R_{max}$ (GtC/year) and $\kappa$ for six different cases. The top plots represent the results of the high-risk scenario ($T_c$ = 2), and the bottom ones show the results of the low-risk scenario ($T_c$ = 3). The plots on the left, middle, and right depict the results for different values of the threshold (d = 1.1, 1.25, and 1.5).}
\label{fig3}
\end{figure}

Social actions can also affect the time to the tipping point. In addition to analysing the evolution of temperature over time for the modified and baseline models, we are also interested in understanding the timing and how it is affected in the modified model under different mitigation strategies. First, we need to define the time to the tipping point. In this work, we measure the time when the temperature anomaly in the modified model becomes \textit{d} times larger than the temperature in the baseline model at the same time, and refer to this specific instant as the \textit{time to the tipping point}. We can now study the evolution of this time as a function of $R_{max}$ and $\kappa$, which cover different scenarios of carbon emissions and mitigation strategies. Figure \ref{fig3} illustrates this analysis. White regions represent cases where the tipping point did not occur. Coloured regions, representing cases with a tipping point, again fit into V-shaped areas similar to the results in \ref{fig2}. Results are provided for two predefined scenarios: high risk (plots A, B, and C) and low risk (plots D, E, and F). By comparing the worst cases of both scenarios, we can see that the tipping point occurs around the year 2120 in the low-risk scenario and around the year 2060 in the high-risk scenario. This suggests that the tipping point can happen 60 years sooner in the high-risk scenario. Additionally, in the low-risk scenario with a higher critical temperature, we observe that the V-shaped regions become smaller, meaning that fewer choices of $R_{max}$ and $\kappa$ lead to a tipping point. We also observe that by selecting different sets of values for $R_{max}$ and $\kappa$, the tipping point can be delayed by around 60 years in the low-risk scenario and 120 years in the high-risk scenario.

\begin{figure}[!t]
\includegraphics[width=\textwidth]{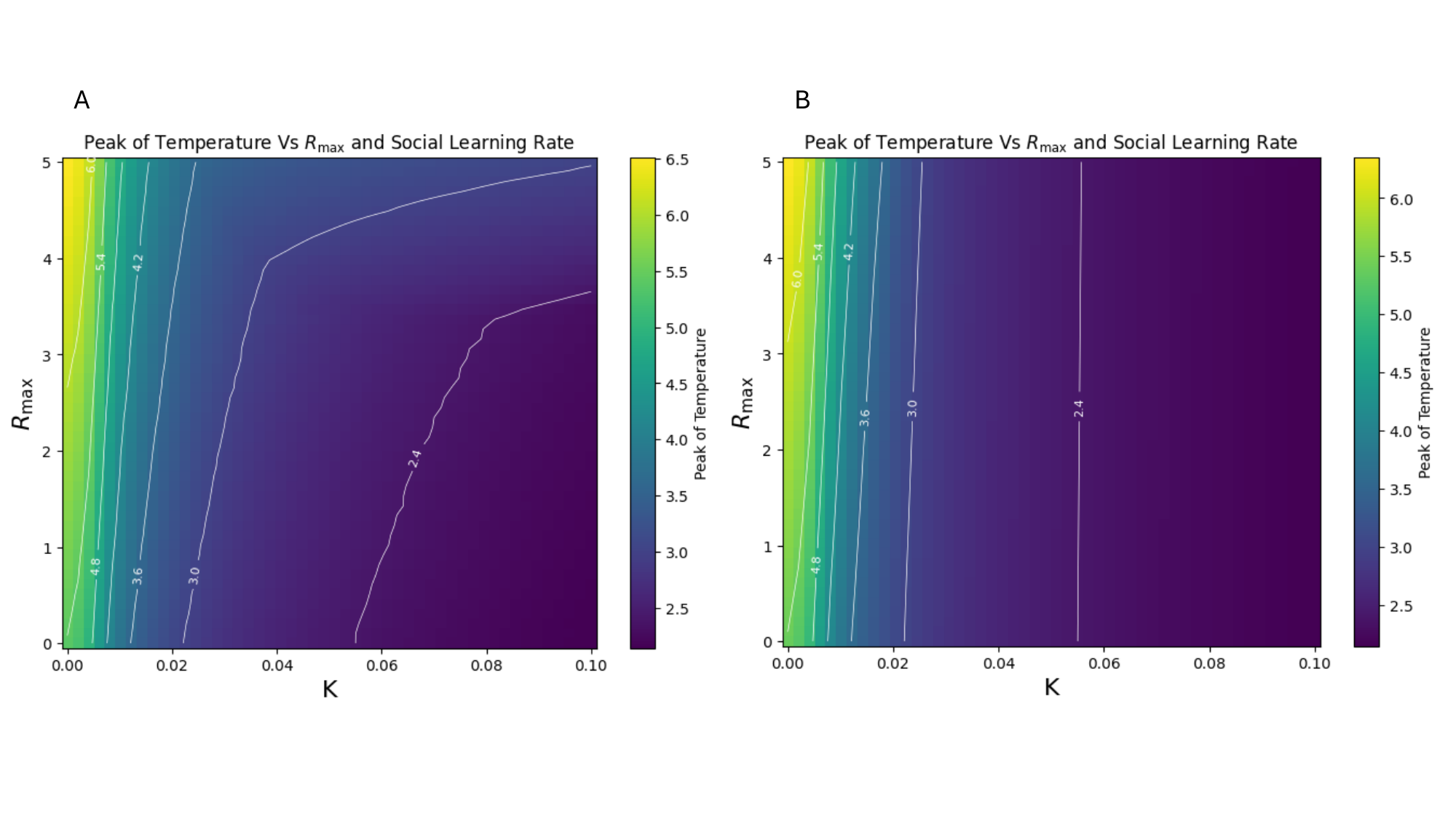}
\caption{\textbf{Peak Temperature Anomaly for Different Emission and Mitigation Scenarios}
Peak temperature anomaly (Celsius) is depicted against $R_{max}$ (GtC/year) and $\kappa$ for the high-risk scenario (left) and the low-risk scenario (right). Isoclines show that as K increases, for the high-risk scenario, the decline in the peak of temperature is lower than low-risk scenario.}
\label{fig4}
\end{figure}

We have also shown that social interventions and the level of carbon emissions can influence the peak temperature deviation. Studying peak temperature deviation is essential as it provides valuable insights into the severity of each scenario and mitigation strategies. The analysis is depicted in Figure \ref{fig4}, where peak temperature anomaly is illustrated against $R_{max}$ and $\kappa$ for two predefined scenarios: high-risk (left) and low-risk (right). Both plots confirm that by increasing $\kappa$ and decreasing $R_{max}$, peak temperature anomaly decreases. Additionally, comparing these two plots reveals that more severe cases tend to occur in high-risk scenarios due to their elevated values for peak temperature anomaly and bent isoclines, suggesting a slower decrease in peak temperature anomaly.


In addition, we have shown that critical temperature plays a pivotal role in determining the severity of our results and distinguishing between scenarios. So far, we have used two values of critical temperature to define high-risk and low-risk scenarios. Returning to the idea of the difference in AUC between the modified and baseline models, we examine the effect of critical temperature, varying it from 1.5 to 5 degrees while $R_{max}$ is set to the constant value of 5 GtC/year. Figure 5 illustrates the outcome of this analysis. The difference in AUC against critical temperature and $\kappa$ is depicted (left), and as shown, it is split into two regions (yellow/green and blue). The yellow/green part represents the cases in which the difference in AUC is high and the tipping term is activated, causing the system to go through a tipping point, while the blue region represents cases where the system does not go through a tipping point. This clearly shows the bifurcated nature of our model. To provide deeper insights, two time series of temperature anomaly versus time are depicted for two examples. The orange circle represents the no-tipping case ($R_{max} = 3$ GtC/year, $\kappa = 0.05$), and the blue star represents the tipping case ($R_{max} = 1.6$ GtC/year, $\kappa = 0.09$). Another conclusion is that for each value of critical temperature, there is a threshold value for $\kappa$, beyond which the modified model does not go through a tipping point. For example, for $T_c = 3$, the tipping point only occurs if $\kappa$ is greater than 0.03.


\begin{figure}[t]
\includegraphics[width=\textwidth]{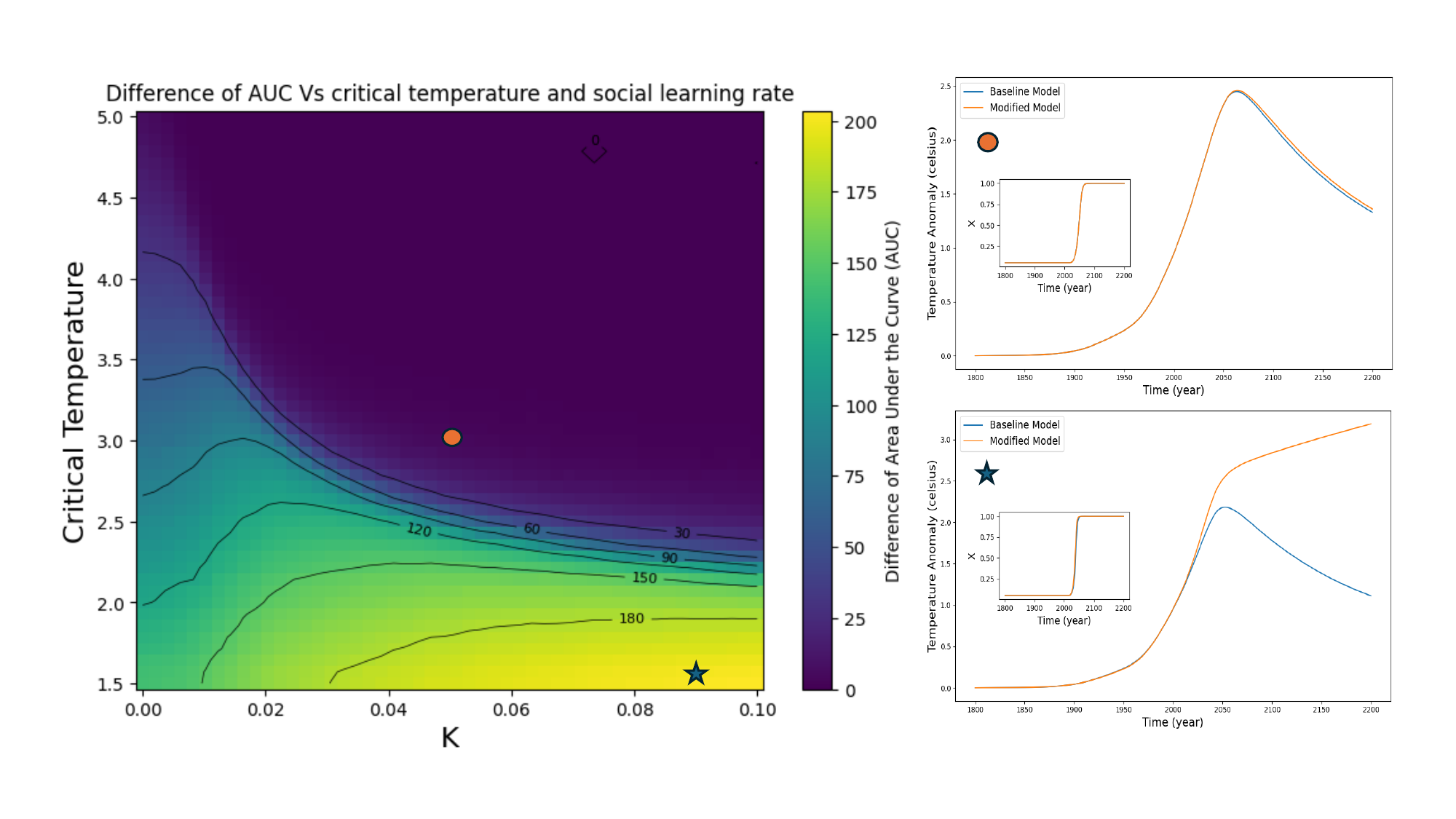}
\caption{\textbf{The difference in AUC versus critical temperature and social learning rate} The Difference in Area under the curve (AUC) is depicted versus critical temperature and $\kappa$ (left). Two examples representing two different choices of $T_c$ and $\kappa$ are chosen and illustrated on RHS. The plot at the top (orange circle) represents a case in which the modified model does not go through a tipping point, and the one at the bottom (blue star) represents a case of tipping point}
\label{fig5}
\end{figure}

In addition to the analysis performed so far, the social model also exhibits a tipping point, and we have shown that the climate tipping point can trigger it. For high delta values, the social model can have up to three equilibrium points (at 0, 1, and one in (0,1)). To simplify the analysis, we define $-\beta + f(T)$ as $\psi$. The number and stability of these equilibrium points vary based on the value of $\psi$.

The social model has stable equilibrium points at 0 and 1, with an unstable point in between for high $\psi$, as shown in Figure 6(A). As $\psi$ decreases, 0 becomes destabilised, leaving the system with a single stable equilibrium point at 1, as depicted in Figure 6(B). Since at the start of the simulations, $x$ is very small and close to zero ($x_0 = 0.05$ in our study), in cases with high $\psi$, the social model converges to the stable equilibrium point at 0 (no mitigation). Conversely, when $\psi$ is low, the social model converges to the only stable equilibrium point at 1 (full mitigation) (full stability analysis of the social model is provided in the Supplementary Materials).

As a result, for specific choices of climate parameters and for high delta ($\delta$ = 3), when the climate model undergoes a tipping point and a sudden increase in temperature, $f(T)$ increases, which raises $\psi$ and destabilises 0, triggering the social model to tip to $x = 1$. The combination of the tipping points in the climate and social models results in an interesting dynamic, as shown in panel C. For high beta (net cost of mitigation), an increase in $f(T)$ during the climate tipping point does not cause $\psi$ to increase sufficiently, and the social model converges to the stable equilibrium point at 0. Consequently, no dynamic interaction is observed. However, for smaller beta, $\psi$ increases sufficiently due to the tipping point in the climate model, causing the social model to shift to the stable equilibrium point at 1, thereby generating diverse dynamics and interactions between the two models.

\begin{figure}[t]
\includegraphics[width=\textwidth]{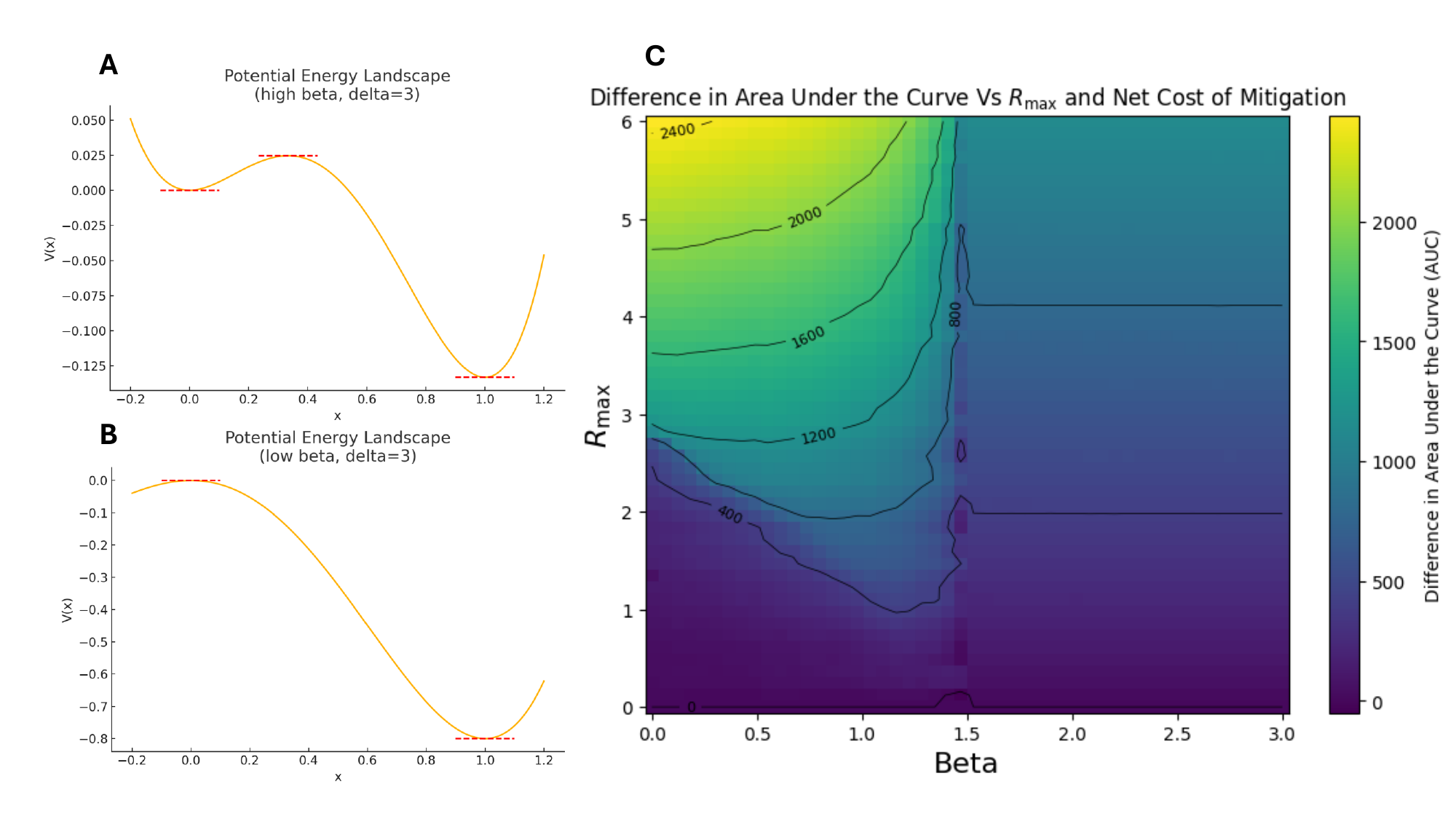}
\caption{\textbf{The Coupled Tipping Dynamics Between Climate and Social Models}
As shown in panels A and B, for high values of delta, decreasing beta destabilises 0 and triggers the social model to tip to x=1. A more detailed illustration of this effect in panel C shows that beyond a specific threshold value of beta, the tipping term does not lead to a significant difference between the modified and baseline model, whereas below that threshold, the tipping element activates.}
\label{fig5}
\end{figure}

\section{Discussion}

In this study, we coupled an Earth System Model (ESM) with a social behaviour model to demonstrate how social actions can delay climate tipping points indefinitely. Our analysis of social parameters revealed that the tipping point can be most effectively averted by accelerating the social learning rate. While other factors, such as the net cost of mitigation and the strength of social norms, influence the climate state, they were found to be largely ineffective in preventing tipping points, but we have also shown that for high values of delta, the tipping point in the climate model can trigger the one in the social model. Additionally, the social learning rate was shown to significantly affect the time to tipping point. For each value of $R_{max}$, representing the strength of the tipping term, we identified a range of social learning rates that can avoid the tipping point and delay its onset. We also illustrated the effect of social learning on the peak of temperature anomalies.

To ensure a thorough analysis, we defined two scenarios based on the critical temperature that triggers the tipping term. It was shown that in high-risk scenarios, the climate state is more severely impacted, and a tipping point is more likely than in low-risk scenarios. Furthermore, we examined the effect of the critical temperature on the likelihood of a tipping point, finding that it plays a crucial role in distinguishing between "tipping" and "no-tipping" cases.

The model and assumptions we used limited the scope of our study. In the social behaviour model, we assumed a homogeneously mixed population, with decisions made individually by comparing the utility of mitigation versus non-mitigation. However, in reality, decisions occur at various societal levels—families, cities, and countries—and people’s choices are not binary but exist on a spectrum. Additionally, the utility of our model is based on the net cost of mitigation, the strength of social norms, and the cost of temperature anomalies. However, other factors, such as advertising and government policies, could also be relevant. Furthermore, human decision-making is not as straightforward as utility comparison; many psychological, social, and cultural factors play crucial roles, and these could be incorporated into future models to enhance realism. Moreover, we aimed to create a simplified model to generate hypotheses and explore new dynamic regimes. Consequently, the Earth-system model and the additional tipping term we employed could be refined in future studies to produce more comprehensive insights and understanding.

In future work, we will focus on addressing the limitations of this study. Specifically, we are interested in exploring the similarities and differences in results by using modified and more complex social behaviour models. Additionally, the climate system can experience various types of bifurcations and tipping points. Therefore, we aim to expand our study to include a broader range of bifurcations and examine their occurrence within an ESM model coupled with a social behaviour model. 

The findings of this study not only enhance our understanding of climate dynamics coupled with social behaviour but also highlight the need to accelerate mitigation efforts and the social learning rate to prevent tipping points. This research demonstrates the crucial role of mitigation in averting climate tipping points. Ultimately, we believe it is vital to study social behaviour dynamics in the context of climate change and tipping points.

\clearpage

\bibliographystyle{plain}

\newpage
\section{Supplementary Materials}

\subsection{Formulation of socio-climate model. }

\subsubsection{behaviour Dynamics}

Individuals in our model are divided into two groups "mitigators" and "non-mitigators". they can move to the other group at the rate of social learning rate. They decide based on the utility of each option. The utility of each option depends on the costs of climate change mitigation, costs associated with the average global temperature anomaly, and the utility of social norms. The utility of being a mitigator is given by

\begin{equation}
    e_M = -\alpha + c \hat{f}(T) + \delta x
\end{equation}

Here, $\alpha$ is the cost of adopting mitigative strategies, c is a proportionality constant, and $\hat{f}(T)$ is the cost of temperature anomaly of T degrees Celsius. $\delta$ is the strength of social norms and x is the fraction of mitigators of the population. The utility of a non-mitigator is given by:

\begin{equation}
    e_N = -\gamma - \hat{f}(T) + \delta (1-x)
\end{equation}

where $\gamma$ is the cost of non-mitigative behaviour. If we show the rate of change between these two groups (social learning rate) with $\kappa$ then the total rate at which non-mitigators switch to being mitigators is:

\begin{equation}
    r_{N->M} = kx(1-x) max(e_M-e_N,0)
\end{equation}

Also, the total rate at which the mitigators switch to being non-mitigators is:

\begin{equation}
    r_{M->N} = kx(1-x) max(e_N-e_M,0)
\end{equation}

As a result, the net rate of change is:

\begin{equation}
    \frac{dx}{dt} = r_{N->M} - r_{M->N} = kx(1-x)(e_M - e_N)
\end{equation}

Using the terms we provided for $e_M$ and $e_N$ we will have:
\begin{equation}
    \frac{dx}{dt} =  kx(1-x)(\gamma - \alpha + (c+1)\hat{f}(T)+ \delta (2x-1))
\end{equation}

which can be rewritten as:

\begin{equation}
    \frac{dx}{dt} =  kx(1-x)(-\beta + f(T)+ \delta (2x-1))
\end{equation}

Here $\beta = \alpha - \gamma$ which defines as net cost of mitigation and $f(T) = (c+1)\hat{f}(T)$.

We also assume that $f(T)$ is a sigmoid function :
\begin{equation}
    \hat{f}(T) = \frac{\hat{f}_{max}}{1+ e^{-\omega(T-T_c)}}
\end{equation}

Here $f_{max}$ is the maximum cost, $\omega$ is the degree of nonlinearity of the sigmoid function and $T_c$ is the critical temperature about which costs are most sensitive to change. we can combine c and $\hat{f}_{max}$:
\begin{equation}
    f_{max} = (c+1)\hat{f}_{max}
\end{equation}

and rewrite the $f(T)$ as:

\begin{equation}
    f(T) = \frac{f_{max}}{1+e^{-\omega(T-T_c)}}
\end{equation}

\subsubsection{Earth system model}

We have used the earth system model \cite{lenton2000land} combined with reduced ocean dynamics \cite{muryshev2015lag}. The full earth system model is as follows:

\begin{align}
    \frac{dC_{at}}{dt} &= \epsilon(t) - P + R_{veg} + R_{so} - F_{oc} \\
    \frac{dC_{oc}}{dt} &= F_{oc} \\
    \frac{dC_{veg}}{dt} &= P - R_{veg} - L \\
    \frac{dC_{so}}{dt} &= L - R_{so} \\
    c\frac{dT}{dt} &= (F_d - \sigma T^4)a_E
\end{align}

P represents photosynthesis which takes the following form:

\begin{equation}
    P(C_{at},T) = k_p C_{v0} k_{MM} (\frac{pCO_{2a}-k_c}{K_M + pCO_{2a}-k_c})(\frac{(15+T)^2 (25-T)}{5625}
\end{equation}

for $pCO_{2a}>= k_c$ and $-15 <= T <= 25$, and zero otherwise. $PCO_{2a}$ is defined as the ratio of moles of $CO_2$ in the atmosphere to the total number of moles of molecules in the atmosphere $k_a$:

\begin{equation}
    pCO_{2a}  = \frac{f_{gtm}(C_{at}+C_{at0})}{k_a}
\end{equation}

where $f_{gtm}=8.3259 * 10^{13}$ is the conversion factor from GtC to moles of carbon and $C_{at0}$ is the initial level of $CO_2$ in the atmosphere. 

Plant respiration takes the form:
\begin{equation}
    R_{veg}(T,C_{veg})  = k_r C_{veg} k_A e^{\frac{-E_a}{R(T+T_0)}}
\end{equation}

Soil respiration takes the form :
\begin{equation}
    R_{so}(T,C_{so})  = k_{sr} C_{so} k_B e^{\frac{-308.56}{T+T_0-227.13)}}
\end{equation}

Turnover (constant fraction of plants dying in a given unit of time) takes the form:
\begin{equation}
    L(C_{veg}) = k_t C_{veg}
\end{equation}

flux of $CO_2$ from the atmosphere to the ocean takes the form

\begin{equation}
    F_{oc}(C_{at},C_{oc}) = F_{0\chi}(C_{at} - \zeta \frac{C_{at0}}{C_{oc0} C_{oc}})
\end{equation}

where $\chi$ is the characteristic solubility of $CO_2$ in water and $\zeta$ is the evasion factor.

The net downward flux of absorbed radiation at the surface is :

\begin{equation}
    F_d = \frac{(1-A)S}{4}(1+ \frac{3\tau}{4})
\end{equation}

where A is the surface albedo, S is the incoming solar flux and $\tau$ is the vertical opacity of the greenhouse atmosphere. The opacity of each greenhouse is given by:

\begin{align}
    \tau(CO_2) &= 1.73(pCO_2)^{0.263} \\
    \tau(H_2O) &= 0.0126 (H P_0 e^{\frac{-L}{RT}})^{0.503} \\
    \tau(CH_4) &= 0.0231
\end{align}

H is the relative humidity, $P_0$ is the water vapor saturation constant, L is the latent heat per mole of water, and R is the molar gas constant.

\subsection{Stability Analysis of the social model}

The social model is given by:
\begin{equation}
    \frac{dx}{dt} =  kx(1-x)(-\beta + f(T)+ \delta (2x-1))
\end{equation}

To simplify the analysis, we define $-\beta + f(T)$ as $\psi$. This model has two equilibrium points at 0 and 1 and also another equilibrium point exists if:
\begin{equation}
    \psi+ \delta (2x-1) = 0
\end{equation}

This equilibrium point happens in the range of (0,1) if 
\begin{equation}
    -\delta < \psi < \delta 
\end{equation}

To evaluate the stability of the equilibrium points, we drive the second derivative of x with respect to time:

\begin{equation}
   k \psi +  k \delta(4x-1) - 2 k \psi x - k \delta (6x^2 - 2x)
\end{equation}

For x=0, this becomes:

\begin{equation}
    k ( \psi - \delta)
\end{equation}

So if $\psi > \delta$ it destabilises x=0.

For x=1, this becomes:

\begin{equation}
    -k ( \psi + \delta)
\end{equation}

if $\psi< - \delta$, x=1 will be unstable.

As a result of a positive and high value of $\delta$, if $\psi < \delta$, then three equilibrium points exist at $x=0$, $x=1$, and one more in $[0,1]$. In this case, $x=0$ and $x=1$ are stable, and the third equilibrium point is unstable. If the climate system goes through a tipping point and a sudden increase in temperature, $f(T)$ will increase, leading to an increase in $\psi$. Due to this increase, if $\psi > \delta$, the third equilibrium point collides with $x=0$ and gets annihilated. Additionally, $x=0$ becomes unstable while $x=1$ remains a stable equilibrium point. As a result, the social model converges to $x=1$ for $\psi > \delta$ and to $x=0$ for $\psi < \delta$, as the initial value of $x$ at the beginning of the simulation is very close to zero.

\subsection{Table. Labels for state variables and dynamic processes in the socio-climate model.}
\noindent\resizebox{\textwidth}{!}{%
\begin{tabular}{||c c c c||} 
 \hline
 Parameter & Definition & Baseline Values / Intervals &  Unit \\ [1ex] 
 \hline\hline
 $\kappa$ & Social Learning Rate & (0, 0.2) & $yr^{-1}$ \\ 
 \hline
 $\delta$ & Strength of Social Norm & (0,2) & 1 \\
 \hline
 $\beta$ & Net Cost of Mitigation & (0,2) & 1 \\
 \hline
 $f_{max}$ & maximum of warming cost function f(T) & (5) & 1 \\
 \hline
 $w$ & nonlinearity of warming cost function f(T) & (3) & $K^{-1}$ \\ 
 \hline

 $T_{lim}$ & limit temperature of f(T) & (1.5) & K \\
  \hline

 $R_{max}$ & maximum of function R(T) & (0,5) & GtC/year \\
  \hline

 $R_{0}$ & nonlinearity of R(T) & (5) & $K^{-1}$ \\
  \hline
 $T_c$ & critical temperature R(T) & (2,3) & K \\
  \hline
 $x_0$ & initial proportion of mitigators   & (0.05) & 1 \\[1ex] 
     \hline
 $T_R$ & Temperature of freezing point for water   & 273.15 & K \\[1ex]
      \hline
 $\sigma$ & Stefan–Boltzmann constant   & $5.67 * 10^{-8}$ & $Wm^-2 K^-4$ \\[1ex]
       \hline

  $C_{at0}$ & initial $CO_2$ in atmosphere   & $596$ & $GtC$ \\[1ex]
       \hline
  $C_{ao0}$ & initial $CO_2$ in ocean reservoir   & $1.5 * 10^5$ & $GtC$ \\[1ex]
       \hline
  $C_{veg0}$ & initial $CO_2$ in vegetation reservoir   & $550$ & $GtC$ \\[1ex]
       \hline 
  $C_{so0}$ & initial $CO_2$ in soil reservoir   & $1500$ & $GtC$ \\[1ex]
       \hline
  ${k_p}$ & photosynthesis rate constant  & $0.184$ & ${yr}^{-1}$ \\[1ex]
       \hline    
  ${k_{MM}}$ & photosynthesis normalising constant  & $1.478$ & $1$ \\[1ex]
       \hline  
  ${k_{c}}$ & photosynthesis compensation point  & $29 * 10^{-6}$ & $1$ \\[1ex]
       \hline  
  ${k_{M}}$ & half-saturation point for photosynthesis & $120 * 10^{-6}$ & $1$ \\[1ex]
       \hline
  ${k_{a}}$ & mole volume of atmosphere & $1.773 * 10^{20}$ & $Moles$ \\[1ex]
       \hline
  ${k_{r}}$ & plant respiration constant & $0.092$ & ${yr}^{-1}$ \\[1ex]
       \hline
  ${k_{A}}$ & plant respiration normalising constant & $8.7039 * 10^9$ & $1$ \\[1ex]
       \hline
  ${E_a}$ & plant respiration activation energy & $54.83$ & $J {mol}^{-1}$ \\[1ex]
       \hline
  ${k_{sr}}$ & soil respiration rate constant & $0.034$ & ${yr}^{-1}$ \\[1ex]
       \hline
  ${k_B}$ & soil respiration normalising constant & $157.072$ & $1$ \\[1ex]
       \hline
  ${k_t}$ & turnover rate constant & $0.092$ & ${yr}^{-1}$ \\[1ex]
       \hline  
  ${c}$ & specific heat capacity of Earth’s surface & $4.69*10^{23}$ & $J K^{-1}$ \\[1ex]
       \hline  
  ${a_E}$ & Earth’s surface area & $5.101*10^{14}$ & $m^2$ \\[1ex]
       \hline 
  ${L}$ & latent heat per mole of water & $43655$ & ${mol}^{-1}$ \\[1ex]
       \hline
  ${R}$ & molar gas constant & $8.314$ & $J {mol}^{-1} K^{-1}$ \\[1ex]
       \hline

  ${H}$ & relative humidity & $0.5915$ & $1$ \\[1ex]
       \hline
  ${A}$ & surface albedo & $0.225$ & ${yr}^-1$ \\[1ex]
       \hline
  ${S}$ & solar flux & $1368$ & $Wm^{-2}$ \\[1ex]
       \hline
  ${\tau(CH_4)}$ & methane opacity & $0.0231$ & $1$ \\[1ex]
       \hline
  ${P_0}$ & water vapor saturation constant & $1.4 * 10^{11}$ & $Pa$ \\[1ex]
       \hline
  ${F_0}$ & ocean flux rate constant & $2.5 * 10^{-2}$ & ${yr}^{-1}$ \\[1ex]
       \hline
  ${\chi}$ & characteristic CO2 solubility & $0.3$ & $1$ \\[1ex]
       \hline
  ${\zeta}$ & evasion factor & $50$ & $1$ \\[1ex]
       \hline
  ${s}$ & half-saturation time for $\epsilon(t)$ from 2014 & $50$ & $yr$ \\[1ex]
       \hline
  ${\epsilon_{max}}$ & maximum change in $\epsilon(t)$ from 2014 & $7$ & $GtC {yr}^{-1}$ \\[1ex]
       \hline
\label{tab: abr table}
\end{tabular}%
}

\subsection{Sensitivity Analysis}
\begin{figure}[t]
\hspace*{-2cm} 
\includegraphics[width=1.2\textwidth]{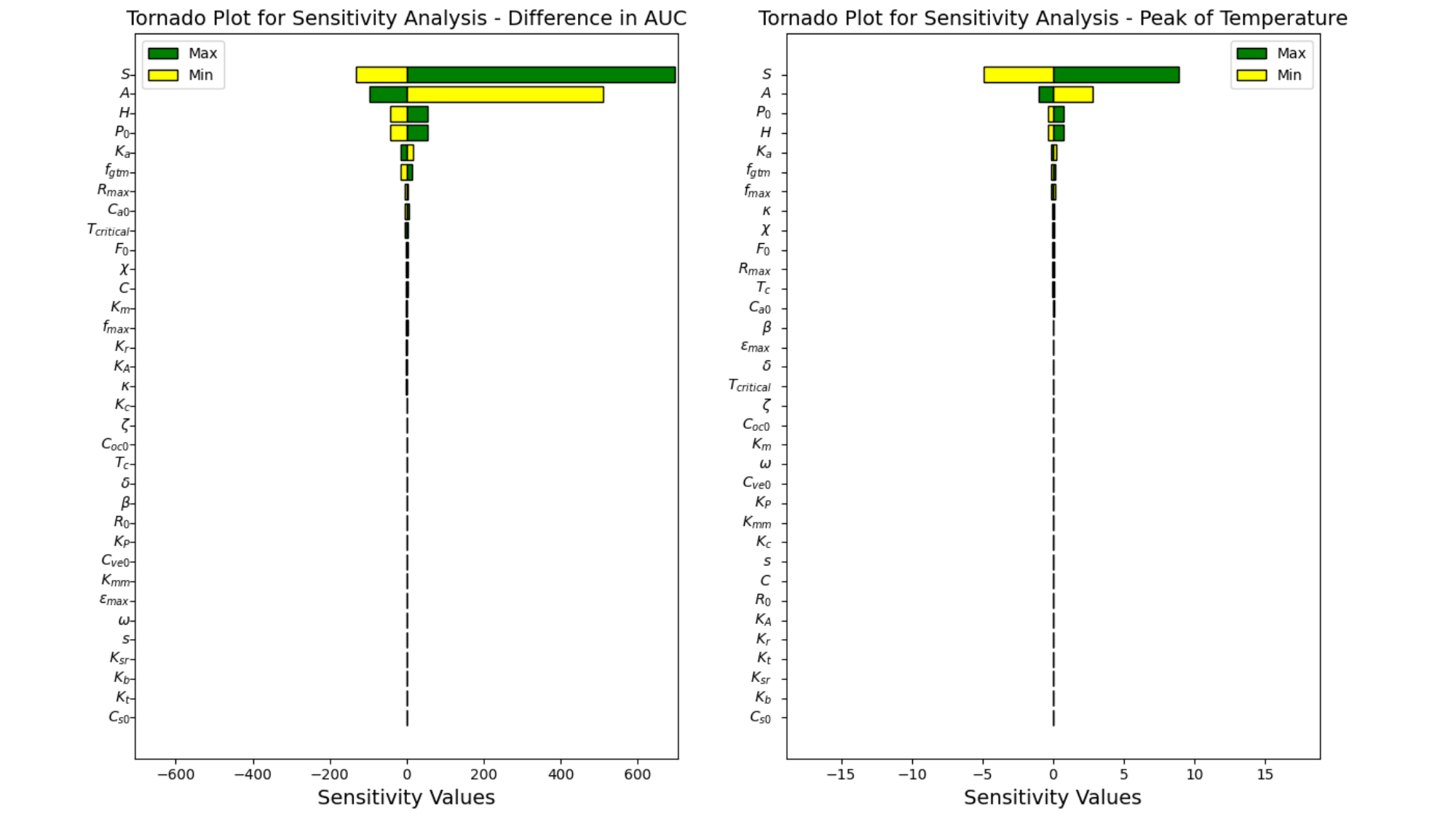} 
\caption{Tornado plot showing the sensitivity analysis for the Difference in the area under the curve (left) and the peak of temperature (right). The yellow (green) bar shows the sensitivity values when the specified parameter has its upper (lower) bound value. The definitions of all parameters and variables are provided in a table in 5.2.}
\label{fig6}
\end{figure}

We have performed a sensitivity analysis to understand the influence of each parameter on the results provided in the previous figures.
Figure \ref{fig6} shows the sensitivity analysis result in a tornado plot for the difference in area under the curve and the peak of temperature anomaly. We have varied all the parameters by 5 percent of their baseline values to generate these plots. By observing both of the plots we can see that Solar flux, initial temperature, Albedo, and Humidity level have the highest influence on the value of the results we have reported. 
This analysis, also reveals the non-monotonic relation between the Difference in AUC and the peak of temperature against the parameters. For instance, an increase in solar flux will increase the difference in AUC while the same result can be achieved by decreasing the Albedo. Additionally, we observe that the sensitivity values are not symmetric, meaning that the constant variation of parameters results in different amounts of change in the difference in AUC and the peak of the temperature anomaly. 

\newpage
\subsection{Historical record of $CO_2$ emission}

\begin{figure}[h]
\centering
\includegraphics[width=0.8\textwidth]{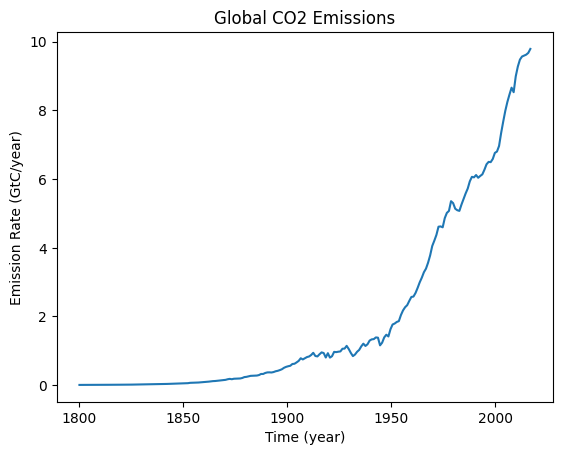} 

\caption{Historical record for Emission rate Vs Time is plotted. This data has been gathered from 1800 until 2017.}
\label{Mu Vs t}
\end{figure}

\subsection{Data and Code Availability}

All data and code supporting this study are available at: 
\newline
https://github.com/Yazdan-Babazadeh/Socio-Climate-model

\end{document}